\documentclass[11pt]{amsart}

\usepackage{amsthm,amsmath,amssymb,amscd,amsfonts}

\theoremstyle{plain}
\newtheorem{theorem}{Theorem}[section]

\newtheorem{proposition}[theorem]{Proposition}

\newtheorem{corollary}[theorem]{Corollary}
\newtheorem{def-thm}[theorem]{Definition-Theorem}
\newtheorem{lemma}[theorem]{Lemma}

\theoremstyle{definition}
\newtheorem{definition}[theorem]{Definition}

\newtheorem{remark}[theorem]{Remark}

\renewcommand{\P}{\mathbb{P}}

\newcommand{\PP}{\mathbb{P}}

\newcommand{\OO}{{\mathcal O}}

\newcommand{\XX}{{\mathcal X}}
\newcommand{\YY}{{\mathcal Y}}

\newcommand{\FF}{{\mathcal F}}
\newcommand{\GG}{{\mathcal G}}

\newcommand{\pr}{{\rm pr}}

\DeclareMathOperator{\Chow}{{Chow}}

\begin{document}

\title[Uniformly effective boundedness]{Uniformly effective boundedness of Shafarevich Conjecture-type for families of canonically polarized manifolds}

\author{Gordon Heier}
\address{Department of Mathematics\\ University of Houston\\ 4800 Calhoun Road, Houston, TX 77204\\USA}

\email{heier@math.uh.edu}

\subjclass[2000]{14C05, 14J10}

\begin{abstract}
The main result of this note is an effective uniform bound for the number of deformation types of certain non-isotrivial families of canonically polarized manifolds. It extends the author's earlier such bound for the classical Shafarevich Conjecture over function fields to the case of higher dimensional fibers, but without the disadvantageous iterated use of Chow or Hilbert varieties that was the core of the proofs in the earlier approaches.
\end{abstract}

\maketitle

\section{Introduction and statement of the effective result}\label{shafsection}
\subsection{The classical case}\label{classical_case}
Let $B$ be a smooth complex projective curve of genus $q\geq 0$. Let
$S\subset B$ be a finite subset of cardinality $s$. The following statement
was conjectured by Shafarevich and proven by Parshin (\cite{P}, case
$S=\emptyset$) and Arakelov (\cite{A}).
\begin{theorem}[\cite{P}, \cite{A}]\label{shaf}
Let $g \geq 2$ be an integer. Then there are only a finite number of isomorphism classes of
nonisotrivial minimal families $f:X\to B$ of curves of genus $g$ with $X$
smooth such that $f:X\backslash f^{-1}(S)\to B\backslash S$ is a smooth
family.
\end{theorem}
Recall that a family of varieties is called {\it isotrivial} if generic fibers are isomorphic to each other.\par
In the article \cite{C}, Caporaso proved that the number of nonisotrivial families in Theorem \ref{shaf} can be bounded by a uniform (but inherently ineffective) constant depending only on $(g,q,s)$. \par
Finally, the article \cite{heier_jmpa} proved an effective uniform bound for the number of nonisotrivial families in Theorem \ref{shaf}, again depending only on $(g,q,s)$. The effectivity is achieved by using Chow varieties instead of stable reduction and $\bar M_g$ as in \cite{C}. For the exact statement, we refer to \cite[Theorem 1.2]{heier_jmpa}. The following is a brief outline of the strategy of proof used in \cite{heier_jmpa}.\par
We assume $q\geq 2$, because the remaining cases $q=0,1$ can simply be handled by executing a degree $2$ base change to a curve of genus $2$. \par
Let $f:X\to B$ be one of the families under discussion. Since $q\geq 2$, $X$
is a manifold of general type and there exist a finite number of rational
$(-2)$-curves $C_i$ $(i\in I)$ with $\bigcup_{i\in I}C_i\subset f^{-1}(S)$
such that, by a well-known theorem of Bombieri, the $5$-canonical map
\begin{equation*}
\varphi_{|5K_X|}:X\to \P(H^0(X,\OO_X(5K_X)))=:\P^{m_X}
\end{equation*}
is an embedding on the complement of $\bigcup_{i\in I}C_i$ and contracts the
$C_i$ to rational double point singularities. In \cite[Proposition 2]{P} it
is stated that $m_X\leq 50\cdot 5^2(gq+s)=1250(gq+s)=:m$ as a consequence of
the Riemann-Roch theorem. We can assume (after linear inclusions) that $m_X=m$ for all families.\par
Since the degree of the divisor $5K_X$ on the smooth fibers $X_b$ is equal to $d:=5(2g-2)$
(independently of $b$), there corresponds to every family $f:X\to B$ a
canonical morphism $\psi_X^0:B\backslash S\to\Chow_{1,d}(\P^m)$ given by
$b\mapsto [\varphi_{|5K_X|}(X_b)]$. We follow the standard convention that $\Chow_{1,d}(\P^m)$ denotes the Chow variety of $1$-dimensional cycles of degree $d$ in $\P^m$. Since $B$ is a smooth curve and $\Chow_{1,d}(\P^m)$ is projective, there exists a unique extension
$\psi_X:B\to\Chow_{1,d}(\P^m)$ that coincides with $\psi_X^0$ on
$B\backslash S$.\par
Clearly, no two nonisomorphic families with the properties described in
Theorem \ref{shaf} can correspond to the same isomorphism class of morphisms $\psi:B\to\Chow_{1,d}(\P^m)$.
Thus, we are left with bounding the cardinality of the set of isomorphism classes of such maps
$\psi_X$, which is achieved by bounding the
degree of the graphs $\Gamma_{\psi_X}\subset B\times\Chow_{1,d}(\P^m)$ and
then applying a rigidity argument to $\Gamma_{\psi_X}$ in conjunction with an
estimate on the complexity of a certain Chow variety of cycles in $B\times\Chow_{1,d}(\P^m)$. This finishes the proof.\par
An undesirable feature of the proof outlined above is that it requires considering Chow varieties twice in an iterated fashion (see the last sentence of the sketch). This is similar to Parshin's original boundedness proof and disadvantageous in more than one way. \par
First of all, the numerical bound obtained in \cite{heier_jmpa} is rather large due to this fact. Secondly, a generalization to the analogous situation with higher dimensional fibers is met with numerous technical difficulties. The purpose of this note is to treat the situation with higher dimensional fibers effectively using a Chow variety only once. While the resulting good effective bound should be interesting in its own right, one can take the point of view that the main achievement of this work is the introduction of a technique that allows for the removal of the second use of a Chow variety.

\subsection{The case of higher dimensional fibers}
If the smooth fibers of $f$ in Theorem \ref{shaf} are taken to be canonically polarized compact complex manifolds instead of curves of genus at least $2$, then there is no analogous finiteness statement. The reason is the following. Take a nonisotrivial family of curves $\tilde f:Y\to B$ as above and a smooth complex projective curve $C$ of genus at least $2$. Let $f$ be the map 
$$X=Y\times C \to B, \quad (y,c)\mapsto \tilde f(y).$$
This family of surfaces continues to be nonisotrivial, but it is not rigid, since a deformation of $C$ gives a deformation of $f$ in the sense described below.\par
Here is the set-up which we will consider in this paper. Let $B$ be a smooth complex projective curve. We assume that its genus $q$ is at least $2$, because the cases $q=0,1$ can again be obtained easily by base change. Let $S\subset B$ be a finite subset of cardinality $s$. These choices are considered arbitrary but fixed for the remainder of the paper. Let $X$ be a projective manifold of dimension $n+1$. A nonisotrivial family $f:X\to B$ such that $f:X\backslash f^{-1}(S)\to B\backslash S$ is a smooth family of canonically polarized compact complex manifolds is called {\it admissible}. Note that the condition of canonical polarization is not usually included in the definition of admissibility, but since all our admissible families will have canonically polarized fibers, we include it nevertheless.\par

\begin{definition}
Let $T,\XX$ be connected quasi-projective varieties. A {\it deformation parametrized by $T$} of the admissible family $f:X\to B$ is a holomorphic map $\FF:\XX\to B\times T$ such that $\FF:\FF^{-1}((B\backslash S)\times \{t_0\}) \to (B\backslash S)\times \{t_0\}$ is isomorphic to $f:X\backslash f^{-1}(S)\to B\backslash S$ for some $t_0\in T$ and $\FF:\FF^{-1}(B\times \{t\})\to B\times \{t\}$ is admissible for every $t\in T$.\par
Two admissible families $f_1:X_1\to B,\ f_2:X_2\to B$ are said to be {\it of the same deformation type} if there exist $T, \XX$ as above and $\FF$, a deformation parametrized by $T$ of $f_1$, such that $\FF:\FF^{-1}((B\backslash S)\times \{t_2\}) \to (B\backslash S)\times \{t_2\}$ is isomorphic to $f_2:X_2\backslash f_2^{-1}(S)\to B\backslash S$ for some $t_2\in T$. 
\end{definition}
As for the history of our problem, recall that Bedulev and Viehweg proved the following in \cite{BedVie} under the assumption that the Minimal Model Conjecture holds. Let $f:X\to B$ be an admissible family such that one (and therefore every) smooth fiber has Hilbert polynomial $h$. Then the number of deformation types of admissible families over $B$ whose smooth fibers also have Hilbert polynomial $h$ is finite. Kov\'acs and Lieblich \cite{KL} then showed that this number can uniformly, but ineffectively, be bounded by a constant depending only on $q$ (the genus of $B$), $s$ (the cardinality of $S$) and $h$. \par
The following Theorem is the main result of this note. The constant $m$ is defined in Definition \ref{m_def}. The square brackets indicate moving self-intersection numbers (see \cite[Definition 11.4.10]{PAGII}). We need to resort to moving self-intersection numbers because we are not making any positivity assumption on the canonical divisor $K_X$ aside from the automatic bigness due to the subadditivity theorem of Kawamata \cite{Kawamata_kod_dim} (see also \cite{kollar} and \cite{viehweg_weak_pos}) and the positivity coming from the nonisotriviality. Moreover, for the remainder of the paper, we choose a fixed base point $b_0\in B\backslash S$ and let $F$ denote the fiber $X_{b_0}$. It is well-known that it follows from the Riemann-Roch theorem that the divisor $L_B=(2q+1)b_0$ on $B$ is very ample. 
\begin{theorem}\label{mthm}
The number of deformation types of admissible families $f:X\to B$ with the numerical value of $(mK_X+\tfrac m 2 (q+1) F)^{[n+1]}+n$ no greater than $p_{\text{max}}$ and the numerical value of $(mK_X+(\tfrac m 2 (q+1)+2q+1)F)^{[n+1]}$ equal to $d$ is no more than the integer in expression \eqref{the_bound} on page \pageref{the_bound}.\end{theorem}
\begin{remark}
When $K_X$ is nef, the effective bound in expression \eqref{the_bound} can be estimated from above in terms of $(K_X^{n+1},K_F^n,n,q)$ as explained in Remark \ref{estim_above_in_terms} and Lemma \ref{def_and_bound_d}. It would then clearly be desirable to bound $K_X^{n+1}$ in terms of $(K_F^n,n,q,s)$. In the case of $1$-dimensional fibers, this is done in \cite[Proposition 1]{P}. However, in the higher dimensional situation, it does not seem to be known how to accomplish this.
\end{remark}

\section{Proof of Theorem \ref{mthm}}
Recall the following effective statement for pluricanonical embeddings from \cite[Corollary 4.2]{heier_doc_math}.
\begin{lemma}\label{effpluri}
If $F$ is a compact complex manifold of complex dimension $\eta$ whose canonical divisor $K_F$ is ample, then $mK_F$ is very ample for any integer $m\geq (e+\frac 1 2)\eta^\frac 7 3+\frac 1 2 \eta^\frac 5 3 + (e+\frac 1 2)\eta^\frac 4 3 + 3\eta+ \frac 1 2 \eta^\frac 2 3+5$, where $e\approx 2.718$ is Euler's constant.
\end{lemma}
\begin{definition}\label{m_def}
From now on, let $m$ be the round-up of the constant from Lemma \ref{effpluri} for $\eta=n$. If $m$ is odd, we increase it by one to make it even.
\end{definition}
The following Lemma will allow us to use Lemma \ref{effpluri} on the fibers of an admissible family in order to get an embedding of $X\backslash f^{-1}(S)$ into projective space with effective control.
\begin{lemma}\label{emb_lemma}
Fix an admissible family $f:X\to B$. Let $s_0,\ldots,s_k$ be a basis for $H^0(B,\OO_B(L_B))$. Let $\tilde \sigma_0,\ldots,\tilde \sigma_{\tilde p}$ be a basis for $H^0(F,\OO_F(mK_F))$. Then there is a finite subset $\tilde S\subset B\backslash S$ such that the sections $\tilde \sigma_0,\ldots,\tilde \sigma_{\tilde p}$ can be extended to elements of $H^0(X,\OO_X(mK_X+\frac m 2 (q+1)F))$ (denoted by the same symbols) in such a way that the sections
\begin{eqnarray*}
&&\tilde \tau_{0,0}=(f^*s_0)\tilde \sigma_0,\ldots,\tilde \tau_{k,0}=(f^*s_k)\tilde \sigma_0,\ldots,\\ &&\tilde \tau_{0,{\tilde p}}=(f^*s_0)\tilde \sigma_{\tilde p},\ldots,\tilde \tau_{k,{\tilde p}}=(f^*s_k)\tilde \sigma_{\tilde p}\in \\ && H^0(X,\OO_X(mK_X+(\tfrac m 2 (q+1)+2q+1)F))
\end{eqnarray*}
yield an embedding of $X\backslash f^{-1}(S\cup\tilde S)$ into $\PP^{(k+1)({\tilde p}+1)-1}$ (with homogeneous coordinates $[X_{0,0},\ldots,X_{k,0},\ldots,X_{0,{\tilde p}},\ldots,X_{k,{\tilde p}}]$).
\end{lemma}
\begin{proof} Since $f$ is nonisotrivial, the locally free sheaf $f_*\OO_X(mK_X)$ is ample on $B$ (see, e.g., \cite[Proposition 4.6]{Viehweg}). Moreover, from Siu's invariance of plurigenera \cite{Siu_invar_pluri}, it follows that $h^0(X_b,\OO_{X_b}(mK_{X_b}))$ is constant for $b\in B\backslash S$. Then \cite[Lemma 7.6]{Kovacs_log_van_thm} yields that $K_X$ is ample with respect to $X\backslash f^{-1}(S)$.\par
We now apply Takayama's extension theorem \cite[Theorem 4.1]{Takayama}, which is in turn based on work of Siu \cite{Siu_inv_plurigen_Grauert}. We refer the reader to \cite[Theorem 4.1]{Takayama} for the precise statement and simply note that Takayama's divisor $L$ is $K_X+qF$ in our situation, and $S$ is $F$. According to the extension theorem, the sections $\tilde \sigma_0,\ldots,\tilde \sigma_{\tilde p}$ extend to elements of $H^0(X,\OO_X(mK_X+\tfrac m 2 (q+1)F))$.\par
Lemma \ref{effpluri} tells us that the extended sections $\tilde \sigma_j$ $(j=0,\ldots,{\tilde p})$ yield an embedding of all smooth fibers $X_b$ where $b$ is not in a finite subset $\tilde S$. Namely, we can take this finite subset $\tilde S$ to be the locus of points $b$ where the fiber $X_b$ is smooth, but $\tilde \sigma_0|_{X_b},\ldots,\tilde \sigma_{\tilde p}|_{X_b}$ are not linearly independent.\par
Since $L_B$ is very ample, and by the properties of the Segre embedding, it is now clear that the $\tilde \tau_{i,j}$ yield an embedding of $X\backslash f^{-1}(S\cup \tilde S)$ into $\PP^{(k+1)({\tilde p}+1)-1}$.
\end{proof}
Lemma \ref{emb_lemma} implies the following Proposition, which is obtained by replacing the extended sections $\tilde \sigma_0,\ldots,\tilde \sigma_{\tilde p}$ with a basis $\sigma_0,\ldots,\sigma_{p}$ for $H^0(X,\OO_X(mK_X+\frac m 2 (q+1)F))$. This replacement allows us to remove the set $\tilde S$ in the statement. The proof consists of applying Takayama's extension theorem to each fiber $X_b$ with $b\in \tilde S$ with $L=K_X+(q+1)F-X_b$. The point of inserting the factor $(q+1)$ is to make sure that $(q+1)F-X_b$ is linearly equivalent to a nef effective divisor on $X$. This particular choice of $L$ guarantees that the extended sections are all in $H^0(X,\OO_X(mK_X+\frac m 2 (q+1)F))$, independently of $b\in \tilde S$.
\begin{proposition}
Fix an admissible family $f:X\to B$. Let $s_0,\ldots,s_k$ be a basis for $H^0(B,\OO_B(L_B))$. Let $\sigma_0,\ldots,\sigma_{p}$ be a basis for $H^0(X,\OO_X(mK_X+\frac m 2 (q+1)F))$. Then the sections
\begin{eqnarray*}
&&\tau_{0,0}=(f^*s_0)\sigma_0,\ldots,\tau_{k,0}=(f^*s_k)\sigma_0,\ldots,\\ &&\tau_{0,{p}}=(f^*s_0)\sigma_{p},\ldots,\tau_{k,{p}}=(f^*s_k)\sigma_{p}\in \\ && H^0(X,\OO_X(mK_X+(\tfrac m 2(q+1)+2q+1)F))
\end{eqnarray*}
yield an embedding of $X\backslash f^{-1}(S)$ into $\PP^{(k+1)({p}+1)-1}$ (with homogeneous coordinates $[X_{0,0},\ldots,X_{k,0},\ldots,X_{0,{p}},\ldots,X_{k,{p}}]$). \end{proposition}
We will denote this embedding by $\varphi_1$. By the Riemann-Roch theorem, we have that
\begin{equation*}
k=h^0(B,\OO_B(L_B))-1= q+1.
\end{equation*}
The following Proposition will give an upper bound for $p=h^0(X,\OO_X(mK_X+\frac m 2 (q+1)F))-1$. The corresponding case of a very ample divisor was treated in \cite[Proposition 1.6]{Heier_Math_Nachr}. Similar arguments already appeared in \cite{Matsusaka} and \cite{Lieberman_Mumford}. The square brackets again indicate moving self-intersection numbers (see \cite[Definition 11.4.10]{PAGII}). 
\begin{proposition}\label{h0_bound}
Let $X$ be an $n$-dimensional compact complex manifold and $L$ a big divisor on $X$ such that the linear series $|L|$ induces a birational map from $X$ to projective space. Then  
\begin{equation*}
h^0(X,L)\leq L^{[n]}+n.
\end{equation*}
\end{proposition}
\begin{proof}
We proceed by induction on the dimension, with the case $n=1$ following immediately from the Riemann-Roch theorem.\par
Let the effective divisor $D$ be a general element of the linear series $|L|$. By \cite[Proposition 18.10]{Harris_first_course}, $D$ is irreducible.
It follows from the standard short exact sequence
\begin{equation*}
0\to\OO_X\to\OO_X(L)\to\OO_D(L)\to 0
\end{equation*}
that
\begin{equation*}
h^0(X,\OO_X(L))\leq h^0(X,\OO_X)+h^0(D,\OO_D(L)).
\end{equation*}
Let $\nu:D'\to D$ be a desingularization of $D$. Simply by pulling back sections, we have
\begin{equation*}
h^0(D,\OO_D(L))\leq h^0(D',\nu^*\OO_D(L)).
\end{equation*}
By applying the induction hypothesis to $D'$, we can furthermore conclude that 
\begin{equation*}
 h^0(D',\nu^*\OO_D(L))\leq (\nu^*(L\cap D))^{[n-1]}+(n-1)=L^{[n]}+(n-1).
\end{equation*}
Since $h^0(X,\OO_X)=1$, this proves the Proposition.
\end{proof}
Due to Proposition \ref{h0_bound}, and possibly after a linear inclusion, we can and do assume from now on that $p$ is actually equal to its largest possible value $p_{\text{max}}$. 
\begin{remark}\label{estim_above_in_terms} When $K_X$ is nef, $p_{\text{max}}$ in Theorem \ref{mthm} can be replaced by an explicit bound in terms of $(K_X^{n+1}, K_F^n,n,q)$ obtained as follows: According to Koll\'ar's effective base point freeness theorem \cite[Theorem 1.1]{kollar_bpf}, $\tilde m K_X$ is base point free for
$$\tilde m := 2(n+3)!(n+3).$$
Now both $\tilde m K_X$ and $\tfrac {\tilde m} 2(q+1) F$ are base point free, and the moving self-intersection number of their sum is just the usual one. Consequently, 
\begin{eqnarray*}
&&(mK_X+\tfrac m 2 (q+1)F)^{[n+1]}+n\\
&\leq & (\tilde m K_X+\tfrac {\tilde m} 2 (q+1)F)^{[n+1]}+n\\
&=& (\tilde m K_X+\tfrac {\tilde m} 2 (q+1)F)^{n+1}+n\\
&=&\tilde m^{n+1} K_X^{n+1}+\tfrac 1 2 {\tilde m}^{n+1}(q+1)(n+1)K_X^{n}.F+n\\
&=&\tilde m^{n+1} K_X^{n+1}+\tfrac 1 2 {\tilde m}^{n+1}(q+1)(n+1)K_F^{n}+n.
\end{eqnarray*}
\end{remark}
The following upper bound on the degree of the closure of the image of $X\backslash f^{-1}(S)$ under $\varphi_1$ is also crucial for our effective argument.
\begin{lemma}\label{def_and_bound_d}
The degree of $\overline{\varphi_1(X\backslash f^{-1}(S))}$, i.e.,\ the closure of the image of $X\backslash f^{-1}(S)$ under $\varphi_1$, is
\begin{equation*}
d:=(mK_X+(\tfrac m 2 (q+1)+2q+1)F)^{[n+1]}.
\end{equation*}
When $K_X$ is nef, we have the following upper bound:
\begin{equation*}
d\leq \tilde m^{n+1} K_X^{n+1}+{\tilde m}^{n}(n+1)\left(\tfrac {\tilde m} 2 (q+1)+2q+1\right)K_F^{n}.
\end{equation*}
\end{lemma}
\begin{proof}
It follows straight from the definition of moving self-intersection numbers that the degree of $\overline{\varphi_1(X\backslash f^{-1}(S))}$ is
\begin{eqnarray*}
&&(mK_X+(\tfrac m 2 (q+1)+2q+1)F)^{[n+1]}.
\end{eqnarray*}
With the same argument concerning effective base point freeness as above, this can be estimated as claimed in the Lemma.
\end{proof}
Let $\varphi_2$ be the embedding of $B$ into $\PP^k$ given by the sections $s_0,\ldots,s_k$ of $\OO_B(L_B)$. There is a commutative diagram 
$$\begin{CD}
X\backslash (\{\sigma_0 = 0\}\cup f^{-1}(S)) @>\iota>>X\backslash f^{-1}(S) @>\varphi_1 >> \PP^{(k+1)(p_{\text{max}}+1)-1}\\
@ V f VV @V f VV @V\pi  VV\\
B\backslash S@>=>>  B\backslash S@>\varphi_2|_{B\backslash S}>> \PP^{k}\\
\end{CD}.
$$
Here, $\iota$ denotes the inclusion map and $\pi: \PP^{(k+1)(p_{\text{max}}+1)-1}\to \PP^{k}$ the rational projection map onto the first $k+1$ homogeneous coordinates, which are $X_{0,0},\ldots X_{k,0}$.\par
\begin{lemma}\label{cd_lemma}
The restricted map 
$$\pi: \varphi_1(X\backslash (\{\sigma_0= 0\}\cup f^{-1}(S)))\to \varphi_2(B\backslash S)$$
is holomorphic and can be extended to a holomorphic map
$$\pi:\varphi_1(X\backslash f^{-1}(S))\to \varphi_2(B\backslash S)$$
such that the diagram
$$\begin{CD}
X\backslash f^{-1}(S) @>\varphi_1 >> \varphi_1(X\backslash f^{-1}(S))\\
@V f VV @V\pi  VV\\
B\backslash S@>\varphi_2>> \varphi_2({B\backslash S})\\
\end{CD}
$$
is a commutative diagram of holomorphic maps. In fact, the diagram is an isomorphism of families.
\end{lemma}
\begin{proof}
The map 
$$\pi: \varphi_1(X\backslash (\{\sigma_0= 0\}\cup f^{-1}(S)))\to \varphi_2(B\backslash S)$$
is clearly holomorphic by construction. Let $x\in X\backslash f^{-1}(S)$ be such that $\sigma_0(x)=0$. Then let $\pi(\varphi_1(x))=\varphi_2(f(x))$. All remaining claims are obvious.
\end{proof}
We now turn to the key argument that makes the second use of a Chow (or Hilbert) variety, which was the core of the proofs in the earlier approaches, unnecessary.\par
The degree of $\varphi_2(B)$ in $\PP^{k}$ is equal to $d_B:=2q+1$. It is well-known (and easy to establish by using general linear projections) that there is a finite set of homogeneous polynomials of degree $d_B$, denoted $\{f_\alpha(X_0,\ldots, X_k)\}_{\alpha}$, such that
$$\varphi_2(B)=\bigcap_{\alpha}\{f_\alpha=0\}.$$ \par
We now lift the $f_\alpha(X_0,\ldots, X_k)$ to homogeneous polynomials $\tilde f_\alpha(X_{0,0},\ldots,X_{k,0},\ldots,X_{0,p_{\text{max}}},\ldots,X_{k,p_{\text{max}}})$ by letting
$$\tilde f_\alpha(X_{0,0},\ldots,X_{k,0},\ldots,X_{0,p_{\text{max}}},\ldots,X_{k,p_{\text{max}}})=f_\alpha(X_{0,0},\ldots, X_{k,0}).$$
The set $\{\tilde f_\alpha\}_\alpha$ defines a subvariety $W$ of $\PP^{(k+1)(p_{\text{max}}+1)-1}$, which is simply $\overline{\pi^{-1}(\varphi_2(B))}$.\par
Let ${\Chow}_{\kappa,\delta}(W)$ be the Chow variety of $\kappa$-dimensional subvarieties of $W$ of degree $\delta$. Let ${\Chow}_{\kappa,\delta}'(W)$ denote the union of those irreducible components of ${\Chow}_{\kappa,\delta}(W)$ whose general points represent irreducible cycles. Since 
$$Z:=\overline{\varphi_1(X\backslash f^{-1}(S))}\subset W,$$
$Z$ corresponds to a point in ${\Chow}_{n+1,d}'(W)$, which we will denote $[Z]$. Recall that the value of $d=\deg(Z)$ was determined in Lemma \ref{def_and_bound_d}.\par
In \cite{Guerra} the following Proposition is proven based on an argument from \cite{Kollarbook}.
\begin{proposition}\label{Chowbound}
Let $W\subset \P^N$ be a projective variety defined by equations of degree no more than $\delta_1$. Then the number of irreducible components of $\Chow'_{\kappa,\delta_2}(W)$ is no more than
\begin{equation*}
{(N+1)\max\{\delta_1,\delta_2\}\choose N}^{(N+1)\left(\delta_2 {\delta_2+\kappa-1\choose \kappa }+{\delta_2+\kappa-1 \choose \kappa-1}\right)}.
\end{equation*} 
\end{proposition}
When we apply Proposition \ref{Chowbound} to our situation, we find that $d_B=\delta_1\leq  \delta_2=d$. Therefore, if we let $N=(k+1)(p_{\text{max}}+1)-1$, the number of irreducible components of ${\Chow}_{n+1,d}'(W)$ is no more than
\begin{equation}\label{the_bound}
{(N+1)d\choose N}^{(N+1)\left(d {d+n\choose n+1 }+{d+n \choose n}\right)}.
\end{equation} 
The proof of Theorem \ref{mthm} will therefore be finished once we have proven the following Proposition.
\begin{proposition}
The number of irreducible components of ${\Chow}_{n+1,d}'(W)$ is an upper bound for the number of deformation types in Theorem \ref{mthm}.
\end{proposition}
\begin{proof}
Let $f:X\to B$ be an admissible family. The same arguments as in the proof of Lemma \ref{cd_lemma} show that there is a commutative diagram of holomorphic maps as follows:
$$\begin{CD}
X\backslash f^{-1}(S) @>\varphi_1 >> \varphi_1(X\backslash f^{-1}(S))@>\text{incl.}>>\overline{\varphi_1(X\backslash f^{-1}(S))}=Z\\
@V f VV @V\pi  VV @V\pi VV\\
B\backslash S@>\varphi_2|_{B\backslash S}>> \varphi_2({B\backslash S})@>\text{incl.}>>\varphi_2(B)\\
\end{CD}.
$$\par
There is a Zariski-open neighborhood $U$ of $[Z]$ in a component of ${\Chow}_{n+1,d}'(W)$ that contains $[Z]$ such that for every point $[V]\in U$, the rational map $\pi:\PP^N\to\PP^k$ induces a holomorphic surjective map $\pi: V\to \varphi_2(B)$. We define $\GG':\YY'\to  \varphi_2(B)\times U$, where $\YY'\subset W\times U$ is the universal family over $U$, to be the natural map induced by $\pi$ and the second projection. \par
Let
$$\YY=\left\{(w,u)\in \YY'|\pi(w)\not\in \varphi_2(S)\right\},$$
and let $\GG:\YY\to\varphi_2(B\backslash S)\times U$ be the restriction of $\GG'$. Due to the generic flatness theorem \cite[6.9.1]{Grothendieck} and the openness of the projection, the map $\pr_2:\YY\to U$, possibly after shrinking $U$ to a smaller Zariski-open subset, is smooth.\par
On the other hand, since $\varphi_2(B\backslash S)$ is a smooth curve, it follows from \cite[Lecture 6, Proposition 6]{Mumford} that $\pr_1\circ \GG$ is flat. Consequently, by \cite[2.1.4, 2.1.7]{Grothendieck} $\GG$ is flat and, possibly after further shrinking of $U$, smooth.\par
We can now finish the proof of the Proposition as follows. Take any two admissible families $f_1,f_2$ that give rise to the same component of ${\Chow}_{n+1,d}'(W)$ and respectively yield the above data (distinguished by an index). Let $T=U_1\cup U_2$. Let 
$$\XX=\left(\YY_1'\cap(W\times U_1)\right)\cup\left(\YY_2'\cap(W\times U_2)\right).$$
Let 
\begin{equation*}
\FF:\XX\to \varphi_2(B)\times T
\end{equation*}
be the obvious map induced by $\GG_1'$ and $\GG_2'$. Then $f_1,f_2$ are of the same deformation type due to $\FF$, which is a deformation parametrized by $T$ that connects $f_1$ and $f_2$.
\end{proof}
\section{The improved bound for the classical case}
For the classical rigid case (with $q\geq 2$) as discussed in Subsection \ref{classical_case}, Theorem \ref{mthm} yields the following Corollary, which significantly improves \cite[Theorem 1.2]{heier_jmpa}. Note that the set of families that is bounded in Corollary \ref{1d_case} is smaller than the set of families in \cite[Theorem 1.2]{heier_jmpa}. However, it is merely a trivial exercise to recover a bound for the same set as in \cite[Theorem 1.2]{heier_jmpa} from Corollary \ref{1d_case} by summing over the possible values of $K_X^2$ ranging from $1$ to $48g(2g+2q+s)$. This upper bound of $K_X^2$ can be found in \cite[Proposition 1]{P}.\par
\begin{corollary}\label{1d_case}
The number of isomorphism classes of admissible minimal families $f:X\to B$ of curves of genus $g\geq 2$ with fixed numerical invariant $K_X^2$ is at most
\begin{equation*}
{(N_1+1)d_1\choose N_1}^{(N_1+1)\left(d_1 {d_1+1\choose 2 }+d_1+1\right)},
\end{equation*}
where
$$d_1=3456g^2+3576gq+1728gs-120q+96g-96,$$
and
$$N_1 =(q+2)(3456g^2+3528gq+1728gs-72q+72g-70)-1.$$
\end{corollary}
\begin{proof}In order to prove the Corollary, we simply have to bound the constants appearing in \eqref{the_bound}. Note that the minimality assumption in the Corollary implies that $K_X$ is nef.\par
Since $3K_F$ is very ample for $1$-dimensional fibers and due to Bombieri's theorem mentioned in the Introduction, there is no need to use the effective base point freeness theorem, and we can simply let $\tilde m=6$. According to Lemma \ref{def_and_bound_d}, $d$ can be estimated from above by 
$$36K_X^{2}+24(5q+4)(g-1)=36K_X^{2}+24(5gq-5q+4g-4).$$
It is known (see \cite[Proposition 1]{P}) that $K_X^2\leq 48g(2g+2q+s),$ so
$$d\leq 3456g^2+3576gq+1728gs-120q+96g-96=:d_1.$$
Moreover, by the Riemann-Roch theorem,
$$k=h^0(B,\OO_B(L_B))-1= q+1,$$
and, by Remark \ref{estim_above_in_terms}, $p_{\text{max}}$ can be replaced by
$$36K_X^{2}+72(q+1)(g-1)+1\leq 3456g^2+3528gq+1728gs-72q+72g-71.$$
Therefore,
\begin{eqnarray*}
N&=&(k+1)(p_{\text{max}}+1)-1\\
&\leq& (q+2)(3456g^2+3528gq+1728gs-72q+72g-70)-1=:N_1.
\end{eqnarray*}
\end{proof}

\end{document}